# Cauchy type problems of linear Riemann-Liouville fractional differential equations with variable coefficients

Myong-Ha Kim, Guk-Chol Ri, Gum-Song Choe, Hyong Chol O*

*Faculty of mathematics,* **Kim Il Sung** *University, Pyongyang DPRK*

*Corresponding author Email:* ryongnam4@yahoo.com

**Abstract:** The existence of solutions to Cauchy type problems of linear Riemann-Liouville fractional differential equations with variable coefficients is considered in a space of integrable functions. First, we consider the existence and uniqueness of solution for a Cauchy type problem with special initial conditions in the space of integrable functions. Then we provide an example of the problem that has no solution in the space of integrable functions. We give a solving method and a representation of solutions for the Cauchy type problem. Last we give some examples.

## 1. Introduction

In recent years, the research on the fractional calculus is becoming an interesting and useful topic among applied scientists.([3,4,8,17,23,24,32,36,37]) Many physical and chemical processes are described by fractional differential equations.([8, 22])

The existence and uniqueness of the solution to an Cauchy type problem for the fractional differential equations were studied in many papers.([1,11,13,16,20,29,35])

In researches for the solving methods and representations of the solutions to Cauchy type problems mainly the linear fractional differential equations with constant coefficients were considered while those with variable coefficients were partly considered for some special forms of equations([2, 3, 9, 10, 14, 15, 16, 25, 29, 30, 33]). In [16, 30, 32, 33] linear fractional differential equations with constant coefficients were considered using Laplace transform and in [6, 7, 9, 16, 21, 29] considered using operational method. Especially, in [16, 32, 34] a representation of Green's function for linear fractional differential operator with constant coefficients in terms of multi-variable Mittag-Leffler function was provided using Laplace transform. In [16, 20] a solution representation of linear Caputo fractional differential equation with constant coefficients was provided using operational calculus of Mikusinski's type. In [9, 29] a solution representation of linear generalized Riemann-Liouville fractional differential equation with constant coefficients was provided using operational calculus of Mikusinski's type. Solving methods and solution representations of linear fractional differential equations with constant coefficients were provided using the method of distribution theory in [27, 28], Neumann series method in [27] and Adomian decomposition method in [10].

In [2] the solving method by Green's function of the system of linear fractional differential equations with constant coefficients was provided and in [30] a representation of Green's function for linear fractional differential operator with variable coefficients was given. In [14], the analytic solution of linear fractional differential operator with variable coefficients was given by power series method. In [5], the analytic solution of a class of fractional differential equations with variable coefficients was provided using properties of Laguerre derivatives and Caputo fractional derivatives.

This article concerns with the existence and representation of solutions to Cauchy type problems of linear Riemann - Liouville fractional differential equations with variable coefficients in a space of integrable functions. First, we consider the existence and uniqueness of solution for a Cauchy type problem with special initial conditions in the space of integrable functions. Then we provide an example of the problem that has no solution in the space of integrable functions. We give a solving method and a representation of solutions for such Cauchy type problems. Last we give some examples which provide counter examples of corollary 3.6 at page 158 of [16].

Our method is just the same method of [16]. We only took notice that some terms are not integrable.

## 2. Preliminaries

Here the definitions and properties of fractional derivatives are described based on [16, 17, 23, 31, 35].

For real numbers $a$, $b$ with $a < b$, we use notations
$$_aR_b := [a, b],\ _aR := [a, \infty),\ R = (-\infty, \infty).$$

We use notations $Z$ and $N$ to denote the sets of all integers and all natural numbers. For $a, b \in Z$ with $a \leq b$, $_aZ_b$ is the set of all integers $i$ satisfying $a \leq i \leq b$ and $_aZ$ the set of all integers $i$ satisfying $i \geq a$. Thus $N = {_1Z}$. Similarly, we use notations $_aN_b$, $_aN$ for $a, b \in N$ with $a \leq b$.

We denote by $C^n_\gamma[a, b]$ the sets of functions satisfying $(x-a)^\gamma f^{(n)}(x) \in C[a, b]$ for $f: (a, b) \to R$ and $0 \leq \gamma < 1$. When $n = 0$, we denote $C_\gamma[a, b] := C^0_\gamma[a, b]$.

Let $\Omega = [a, b]\ (-\infty \leq a < b \leq \infty)$ be a finite or infinite interval of the real axis $R$. We denote by $L_p(a, b)$ the set of those Lebesgue complex-valued measurable functions $f$ for which $\|f\|_p < \infty$, where

$$\|f\|_p := \left(\int_a^b |f(x)|^p dx\right)^{1/p} \quad (1 \leq p < \infty) \tag{2.1}$$

$$\|f\|_p := \operatorname*{ess\,sup}_{a \leq x \leq b} |f(x)| \quad (p = \infty) \tag{2.2}$$

Here $\operatorname*{ess\,sup}_{a \leq x \leq b} |f(x)|$ is the essential maximum of the function $|f(x)|$.

When $p = 1$, we denote $L(a, b) := L_1(a, b)$.

Let $[a, b]\ (-\infty < a < b < \infty)$ be a finite interval and let $AC[a, b]$ be the space of absolutely continuous functions $f$ on $[a, b]$. It is known that $AC[a, b]$ coincides with the space of primitives of Lebesgue summable functions, that is,

$$f(x) \in AC[a, b] \Leftrightarrow f(x) = c + \int_a^x \varphi(t)dt\ (\varphi(x) \in L(a, b)). \tag{2.3}$$

Therefore, an absolutely continuous function $f(x)$ has a summable derivative $f'(x) = \varphi(x)$ almost everywhere on $[a, b]$ and $c = f(a)$.

For $n \in N$ we denote by $AC^n[a, b]$ the space of complex-valued functions $f(x)$ which have all continuous derivatives of order $i$ $(0 \leq i \leq n-1)$ on $[a, b]$ such that $f^{(n-1)}(x) \in AC[a, b]$:

$$AC^n[a, b] := \left\{ f : [a, b] \to C, (D^{n-1}f)(x) \in AC[a, b], \left(D = \frac{d}{dx}\right) \right\}$$

$$AC[a, b] := AC^1[a, b] \tag{2.5}$$

This space is characterized by the following assertion

$$AC^n[a, b] \ni f(x) \Leftrightarrow f(x) = (I_{a+}^n \varphi)(x) + \sum_{k=0}^{n-1} c_k (x-a)^k . \tag{2.6}$$

Here $\varphi(x) \in L(a, b)$ and $c_k$ $(k = 0, 1, \cdots, n-1)$ are arbitrary constants, and

$$(I_{a+}^n \varphi)(x) = \frac{1}{(n-1)!} \int_a^x (x-t)^{n-1} \varphi(t) dt . \tag{2.7}$$

From (2.7) we have

$$\varphi(x) = f^{(n)}(x), \quad c_k = \frac{f^{(k)}(a)}{k!} \quad (k = 0, 1, \cdots, n-1) . \tag{2.8}$$

Let $C$ be the set of all complex numbers. The Riemann-Liouville fractional integrals $I_{a+}^\alpha f$ of order $\alpha \in C$ $(\mathrm{Re}\,\alpha > 0)$ are defined by

$$(I_{a+}^\alpha f)(x) := \frac{1}{\Gamma(\alpha)} \int_a^x \frac{f(t)dt}{(x-t)^{1-\alpha}} \quad (x > a, \; \mathrm{Re}\,\alpha > 0) . \tag{2.9}$$

Here $\Gamma(\alpha)$ is the Gamma function and $(x-t)^{1-\alpha} = e^{(1-\alpha)\ln(x-t)}$. These integrals are called the left-sided fractional integrals. When $\alpha = n \in N$, the definition (2.9) coincides with the nth iterated integral

$$(I_{a+}^n f)(x) = \int_a^x dt_1 \int_a^{t_1} dt_2 \cdots \int_a^{t_{n-1}} f(t_n) dt_n = \frac{1}{(n-1)!} \int_a^x (x-t)^{n-1} f(t) dt \quad (n \in N) . \tag{2.10}$$

The Riemann-Liouville fractional derivatives $D_{a+}^\alpha y$ of order $\alpha \in C$ $(\mathrm{Re}\,\alpha \geq 0)$ are defined by

$$(D_{a+}^\alpha y)(x) := \left(\frac{d}{dx}\right)^n (I_{a+}^\alpha y)(x) = \frac{1}{\Gamma(n-\alpha)} \left(\frac{d}{dx}\right)^n \int_a^x \frac{y(t)dt}{(x-t)^{\alpha-n+1}} \tag{2.11}$$

$$(n = [\mathrm{Re}\,\alpha] + 1, \; x > a).$$

Here $[\mathrm{Re}\,\alpha]$ means the integral part of $\mathrm{Re}\,\alpha$. In particular, when $\alpha = n \in_0 Z$, then

$$(D_{a+}^0 y)(x), \; (D_{a+}^n y)(x) = y^{(n)}(x) . \tag{2.12}$$

Here $y^{(n)}(x)$ is the usual derivative of $y(x)$ of order $n$. If $0 < \mathrm{Re}\,\alpha < 1$, then

$$(D_{a+}^\alpha y)(x) = \frac{1}{\Gamma(1-\alpha)} \frac{d}{dx} \int_a^x \frac{y(t)dt}{(x-t)^{\alpha-[\mathrm{Re}\,\alpha]}} \quad (0 < \mathrm{Re}\,\alpha < 1, \; x > a) \tag{2.13}$$

When $\alpha \in R^+ = (0, \infty)$, then (2.11) take the following form:

$$(D_{a+}^\alpha y)(x) = \frac{1}{\Gamma(n-\alpha)} \left(\frac{d}{dx}\right)^n \int_a^x \frac{y(t)dt}{(x-t)^{\alpha-n+1}} \quad (n = [\alpha]+1, \; x > a) \tag{2.14}$$

and (2.13) is given by

$$(D_{a+}^{\alpha} y)(x) = \frac{1}{\Gamma(1-\alpha)} \frac{d}{dx} \int_a^x \frac{y(t)dt}{(x-t)^{\alpha-[\alpha]}} . \tag{2.15}$$

If $\operatorname{Re}\alpha = 0$ ($\alpha \neq 0$), then (2.11) yields fractional derivatives of a purely imaginary order:

$$(D_{a+}^{i\theta} y)(x) = \frac{1}{\Gamma(1-i\theta)} \frac{d}{dx} \int_a^x \frac{y(t)dt}{(x-t)^{i\theta}} \quad (\theta \in R \setminus \{0\},\ x > a) \tag{2.16}$$

For $\alpha \in C,\ \operatorname{Re}\alpha > 0,\ 1 \le p < \infty$, the spaces of functions $I_{a+}^{\alpha}(L_p)$ is defined by

$$I_{a+}^{\alpha}(L_p) := \{f:\ f = I_{a+}^{\alpha}\varphi,\ \varphi \in L_p(a,\ b)\} . \tag{2.17}$$

For $\alpha \in C\ (\operatorname{Re}\alpha > 0)$, the spaces of functions $L^{\alpha}(a,\ b)$ is defined by

$$L^{\alpha}(a,\ b) := \{y \in L(a,\ b):\ D_{a+}^{\alpha} y \in L(a,\ b)\} . \tag{2.18}$$

**Lemma 2.1** If $\beta \in C\ (\operatorname{Re}\beta > 0)$, then

$$(I_{a+}^{\alpha}(t-a)^{\beta-1})(x) = \frac{\Gamma(\beta)}{\Gamma(\beta+\alpha)}(x-a)^{\beta+\alpha-1} \quad (\operatorname{Re}\alpha > 0) \tag{2.19}$$

$$(D_{a+}^{\alpha}(t-a)^{\beta-1})(x) = \frac{\Gamma(\beta)}{\Gamma(\beta-\alpha)}(x-a)^{\beta-\alpha-1} \quad (\operatorname{Re}\alpha \ge 0) \tag{2.20}$$

In particular, if $\beta = 1$ and $\operatorname{Re}\alpha \ge 0$, then Riemann-Liouville fractional derivatives of constant functions are not equal to zero and we have

$$(D_{a+}^{\alpha} 1)(x) = \frac{(x-a)^{-\alpha}}{\Gamma(1-\alpha)} . \tag{2.21}$$

On the other hand, for $j = 1,\ 2,\cdots,\ [\operatorname{Re}\alpha]+1$, we have

$$(D_{a+}^{\alpha}(t-a)^{\alpha-j})(x) = 0 . \tag{2.22}$$

From (2.22) we have the following result.

**Lemma 2.2** Let $\operatorname{Re}\alpha > 0,\ n = [\operatorname{Re}\alpha]+1,\ n-1 < \operatorname{Re}\alpha < n$. The equality

$$(D_{a+}^{\alpha} y)(x) = 0 \tag{2.23}$$

is valid if and only if

$$y(x) = \sum_{j=1}^{n} c_j (x-a)^{\alpha-j} . \tag{2.24}$$

Here $c_j \in R\ (j = 1,\cdots,\ n)$ are arbitrary constants. In particular, when $0 < \operatorname{Re}\alpha < 1$, the relation $(D_{a+}^{\alpha} y)(x) = 0$ holds if and only if $y(x) = c(x-a)^{\alpha-1},\ \forall c \in R$.

**Lemma 2.3** If $\operatorname{Re}\alpha \ge 0,\ n = [\operatorname{Re}\alpha]+1,\ y(x) \in AC^n[a,\ b]$, then the fractional derivatives $D_{a+}^{\alpha} y$ exist almost everywhere on $[a,\ b]$ and we have

$$(D_{a+}^{\alpha} y)(x) = \sum_{k=0}^{n-1} \frac{y^{(k)}(a)}{\Gamma(1+k-\alpha)}(x-a)^{k-\alpha} + \frac{1}{\Gamma(n-\alpha)} \int_a^x \frac{y^{(n)}(t)dt}{(x-t)^{\alpha-n+1}} . \tag{2.25}$$

**Lemma 2.4.** If $\operatorname{Re}\alpha > 0,\ \operatorname{Re}\beta > 0$ and $f(x) \in L_p(a,\ b)\ (1 \le p \le \infty)$, then we have

$$(I_{a+}^{\alpha} I_{a+}^{\beta} f)(x) = (I_{a+}^{\alpha+\beta} f)(x) \tag{2.26}$$

at almost every point $x \in [a, b]$. If $\alpha + \beta > 1$, then (2.26) is true for any point of $[a, b]$.

**Lemma 2.5.** If $\operatorname{Re}\alpha > 0$, $f(x) \in L_p(a, b)$ $(1 \leq p < \infty)$, then for almost all $x \in [a, b]$ we have

$$(D_{a+}^{\alpha} I_{a+}^{\alpha} f)(x) = f(x). \tag{2.27}$$

From lemma 2.3, lemma 2.4 and lemma 2.5, we have the following lemma.

**Lemma 2.6.** If $\operatorname{Re}\alpha > \operatorname{Re}\beta > 0$, $f(x) \in L_p(a, b)$ $(1 \leq p \leq \infty)$, then we have

$$(D_{a+}^{\beta} I_{a+}^{\alpha} f)(x) = (I_{a+}^{\alpha-\beta} f)(x), \ a.e.\ x \in [a, b]. \tag{2.28}$$

In particular, if $\beta = k \in N$ and $\operatorname{Re}\alpha > k$, then we have

$$(D^k I_{a+}^{\alpha} f)(x) = (I_{a+}^{\alpha-k} f)(x), \ a.e.\ x \in [a, b]. \tag{2.29}$$

**Lemma 2.7** Let $\operatorname{Re}\alpha \geq 0$, $m \in N$, $D = \dfrac{d}{dx}$. If the fractional derivatives $(D_{a+}^{\alpha} y)(x)$, $(D_{a+}^{\alpha+m} y)(x)$ exist, then we have

$$(D^m D_{a+}^{\alpha} y)(x) = (D_{a+}^{\alpha+m} y)(x), \ a.\ e.\ x \in [a, b]. \tag{2.30}$$

**Lemma 2.8.** Let $\operatorname{Re}\alpha > 0$, $n = [\operatorname{Re}\alpha] + 1$, $f_{n-\alpha}(x) := (I_{a+}^{n-\alpha} f)(x)$.

(a) If $1 \leq p \leq \infty$ and $f(x) \in I_{a+}^{\alpha}(L_p)$, then we have

$$(I_{a+}^{\alpha} D_{a+}^{\alpha} f)(x) = f(x). \tag{2.31}$$

(b) If $f(x) \in L(a, b)$, $f_{n-\alpha} \in AC^n[a, b]$, $n-1 < \operatorname{Re}\alpha < n$, then we have

$$(I_{a+}^{\alpha} D_{a+}^{\alpha} f)(x) = f(x) - \sum_{j=1}^{n} \frac{f_{n-\alpha}^{(n-j)}(a)}{\Gamma(\alpha - j + 1)} (x-a)^{\alpha-j}, \ a.e.\ x \in [a, b]. \tag{2.32}$$

If $\alpha = n \in N$, then we have

$$(I_{a+}^{n} D_{a+}^{n} f)(x) = f(x) - \sum_{k=0}^{n} \frac{f^{(k)}(a)}{k!} (x-a)^k, \ a.\ e.\ x \in [a, b]. \tag{2.33}$$

**Remark 2.1** In what follows, we will regard that all the equalities hold almost everywhere on the interval.

## 3. Cauchy type problems of linear equations

For a complex number $\alpha \in C$ with $\operatorname{Re}\alpha \geq 0$, we define the natural number $n$ by

$$n := \begin{cases} [\operatorname{Re}\alpha] + 1, & \alpha \notin N \\ \alpha, & \alpha \in N \end{cases} \tag{3.1}$$

Let $l \in N$ and assume that complex numbers $\alpha_j \in C$ $(j = 1, \cdots, l)$ satisfy

$$0 = \alpha_0 < \operatorname{Re}\alpha_1 < \cdots < \operatorname{Re}\alpha_l < \operatorname{Re}\alpha. \tag{3.2}$$

We consider the following Cauchy type problem of the linear fractional differential equation of order $\alpha \in C$ with variable coefficients;

$$(D_{a+}^{\alpha} y)(x) + \sum_{j=0}^{l} a_j(x)(D_{a+}^{\alpha_j} y)(x) = g(x) \tag{3.3}$$

$$(D_{a+}^{\alpha-k} y)(a+) = b_k, \ k = 1, \cdots, n, \tag{3.4}$$

where the functions $a_j(x)$ and $g(x)$ are complex valued functions with real variable $t$

satisfying

$$a_j(x) \in C[a,\ b]\ (j=0,\cdots,\ l),\ g(x) \in L(a,\ b) \tag{3.4-1}$$

and

$$b_k \in C,\ k=1,\cdots,\ n. \tag{3.4-2}$$

We use the natural number $k_0$ defined by

$$k_0 = \max\left\{k \in_1 N_n : \forall j \in_0 Z_l,\ a_j(t)\frac{(x-a)^{\alpha-\alpha_j-k}}{\Gamma(\alpha-\alpha_j-k+1)} \in L(a,\ b)\right\}. \tag{3.5}$$

### 3.1 A special Cauchy type problem.

We consider a special Cauchy type problem finding the function $y(x)$ satisfying the equation (3.3) and the initial condition

$$(D_{a+}^{\alpha-k}y)(a+) = b_k H(k_0 - k),\ k=1,\cdots,\ n. \tag{3.6}$$

Here

$$H(k) = \begin{cases} 1,\ k \geq 0, \\ 0,\ k < 0. \end{cases} \tag{3.7}$$

We consider the following Volterra integral equation of the second kind corresponding to the Cauchy type problem (3.3), (3.6)

$$\Phi(x) = g(x) - \sum_{j=0}^{l}\sum_{k=1}^{n} b_k H(k_0-k) a_j(x)\frac{(x-a)^{\alpha-\alpha_j-k}}{\Gamma(\alpha-\alpha_j-k+1)} - \sum_{j=0}^{l} a_j(x)(I_{a+}^{\alpha-\alpha_j}\Phi)(x). \tag{3.8}$$

First, we consider the equivalence of the Cauchy type problem (3.3) and (3.6) with the Volterra integral equation (3.8).

**Theorem 3. 1** Assume that for $\alpha, \alpha_j\ (j=0,\cdots,\ l), n,\ a_j(x)(j=0,\cdots,l),\ g(x), b_k\ (k=1,\cdots,n)$ and $k_0$ (3.1), (3.2), (3.4-1), (3.4-2) and (3.5) hold. Then $y(x) \in L^\alpha(a,\ b)$ is the solution to the Cauchy type problem (3.3) and (3.6) if and only if

$$(D_{a+}^\alpha y)(x) = \Phi(x) \in L(a,\ b) \tag{3.9}$$

satisfies the integral equation (3.8).

**Proof.** First we prove the necessity. Let $y(x) \in L^\alpha(a,\ b)$ be a solution to the problem (3.3), (3.6). Then we have

$$(D_{a+}^\alpha y)(x) = g(x) - \sum_{j=0}^{l} a_j(x)(D_{a+}^{\alpha_j}y)(x) \tag{3.10}$$

$$(D_{a+}^{\alpha-k}y)(a+) = b_k H(k_0-k)\ (k=1,\cdots,\ n) \tag{3.11}$$

and $(D_{a+}^\alpha y)(x) \in L(a,\ b),\ g(x) \in L(a,\ b)$. From (3.10) we have

$$\sum_{j=0}^{l} a_j(x)(D_{a+}^{\alpha_j}y)(x) \in L(a,\ b). \tag{3.12}$$

$y(x) \in L^\alpha(a,\ b)$ and (2.18) yield $y(x) \in L(a,\ b)$. By (3.9) and the definition of the fractional derivative, we have

$$(D_{a+}^\alpha y)(x) = (D^n I_{a+}^{n-\alpha}y)(x) = \Phi(x) \in L(a,\ b) \tag{3.13}$$

and from (2.6) we have $y_{n-\alpha}(x) = (I_{a+}^{n-\alpha}y)(x) \in AC^n[a,\ b]$.

The fact that $n-1 < \mathrm{Re}\,\alpha < n$ and lemma 2.8 (a) give

$$(I_{a+}^{\alpha} D_{a+}^{\alpha} y)(x) = y(x) - \sum_{k=1}^{n} y_{n-\alpha}^{(n-k)}(a+) \frac{(x-a)^{\alpha-k}}{\Gamma(\alpha-k+1)}. \qquad (3.14).$$

Since $n-k-1 < \mathrm{Re}(\alpha-k) = \mathrm{Re}\,\alpha - k \le n-k$ we have

$$y_{n-\alpha}^{(n-k)}(x) = (D^{n-k} I_{a+}^{n-\alpha} y)(x) = (D^{n-k} I_{a+}^{n-k-(\alpha-k)} y)(x) = (D_{a+}^{\alpha-k} y)(x).$$

Therefore we have

$$(I_{a+}^{\alpha} D_{a+}^{\alpha} y)(x) = y(x) - \sum_{k=1}^{n} (D_{a+}^{\alpha-k} y)(a+) \frac{(x-a)^{\alpha-k}}{\Gamma(\alpha-k+1)}. \qquad (3.15).$$

Applying the fractional integral operator $I_{a+}^{\alpha}$ in the both sides of (3.9) under consideration of (3.11) and (3.15), we have

$$y(x) = \sum_{k=1}^{n} b_k H(k_0 - k) \frac{(x-a)^{\alpha-k}}{\Gamma(\alpha-k+1)} + (I_{a+}^{\alpha} \Phi)(x) \text{ in } L(a, b). \qquad (3.16)$$

Thus from (3.9) and (3.11), we obtain the relation (3.16). Conversely, we can easily obtain the relations (3.9) and (3.11) from (3.16). Thus the relations (3.9), (3.11) and the relation (3.16) are equivalent. From (3.16), considering (2.20) and (2.28), we have

$$(D_{a+}^{\alpha_j} y)(x) = \sum_{k=1}^{n} b_k H(k_0 - k) \frac{(x-a)^{\alpha-\alpha_j-k}}{\Gamma(\alpha-\alpha_j-k+1)} + (I_{a+}^{\alpha-\alpha_j} \Phi)(x). \qquad (3.17)$$

Therefore, we have

$$\sum_{j=0}^{l} a_j(x)(D_{a+}^{\alpha_j} y)(x) = \sum_{j=0}^{l} \sum_{k=1}^{n} b_k H(k_0 - k) \frac{(x-a)^{\alpha-\alpha_j-k}}{\Gamma(\alpha-\alpha_j-k+1)} + \sum_{j=0}^{l} a_j(x)(I_{a+}^{\alpha-\alpha_j} \Phi)(x). \qquad (3.18)$$

Here the first term of the right side of (3.18) is integrable by the definition (3.5) of $k_0$ and the left side is integrable by (3.12). The second term of the right side is integrable by the assumptions on $a_j$ and $\Phi$. Therefore the equation (3.18) holds in $L(a, b)$.

Substituting (3.9) and (3.18) into (3.10), we have (3.8), i.e. $\Phi(x) \in L(a, b)$ in (3.9) satisfies the integral equation (3.8).

  Now we prove the sufficiency. Assume that $\Phi(x) \in L(a, b)$ satisfies the equation (3.8). The equation (3.8) can be rewritten as

$$\Phi(x) = g(x) - \sum_{j=0}^{l} a_j(x) \left[ \sum_{k=1}^{n} b_k H(k_0 - k) \frac{(x-a)^{\alpha-\alpha_j-k}}{\Gamma(\alpha-\alpha_j-k+1)} - (I_{a+}^{\alpha-\alpha_j} \Phi)(x) \right]. \qquad (3.19)$$

For $\Phi(x) \in L(a, b)$ satisfying (3.8), the equality (3.16) and the equalities (3.9) and (3.11) are equivalent, and thus we have (3.17). Therefore for $\Phi(x) \in L(a, b)$ satisfying (3.8) we have

$$y(x) = \sum_{k=1}^{n} b_k H(k_0 - k) \frac{(x-a)^{\alpha-k}}{\Gamma(\alpha-k+1)} + (I_{a+}^{\alpha} \Phi)(x) \qquad (3.20)$$

$$D_{a+}^{\alpha} y(x) = \Phi(x) \qquad (3.21)$$

$$(D_{a+}^{\alpha-k} y)(a+) = b_k H(k_0 - k) \qquad (3.22)$$

$$(D_{a+}^{\alpha_j} y)(x) = \sum_{k=1}^{n} b_k H(k_0 - k) \frac{(x-a)^{\alpha-\alpha_j-k}}{\Gamma(\alpha-\alpha_j-k+1)} + (I_{a+}^{\alpha-\alpha_j} \Phi)(x) \qquad (3.23)$$

$$y(x) \in L^\alpha(a, b) . \tag{3.24}$$

Substituting (3.21) and (3.23) into (3.19), then we have (3.3) and (3.6).

Thus if $\Phi(x) \in L(a, b)$ satisfies (3.8), then $y(x)$ given by (3.16) satisfies $y(x) \in L^\alpha(a, b)$ and (3.3) and (3.6). Thus the sufficiency is proved, which completes the proof of theorem

In Theorem 3.1 we proved that the Cauchy type problem (3.3), (3.6) and the integral equation (3.8) are equivalent. Therefore in order to prove the existence and uniqueness of solution to the initial vale problem (3.3), (3.6), is sufficient to prove the existence of a unique solution the integral equation (3.8).

**Theorem 3.2.** Under the conditions of Theorem 3.1, there exists a unique solution $y(x) \in L^\alpha(a, b)$ to the Cauchy type problem (3.3), (3.6).

**Proof.** By Theorem 3.1, it is sufficient to establish the existence of a unique solution $\Phi(x) \in L(a, b)$ to the integral equation

$$\Phi(x) = g(x) - \sum_{j=0}^{l}\sum_{k=1}^{n} b_k H(k_0 - k) a_j(x) \frac{(x-a)^{\alpha-\alpha_j-k}}{\Gamma(\alpha-\alpha_j-k+1)} - \sum_{j=0}^{l} a_j(x)(I_{a+}^{\alpha-\alpha_j}\Phi)(x) .$$

Define

$$\Phi_0(x) := g(x) - \sum_{j=0}^{l}\sum_{k=1}^{n} b_k H(k_0 - k) a_j(x) \frac{(x-a)^{\alpha-\alpha_j-k}}{\Gamma(\alpha-\alpha_j-k+1)} , \tag{3.25}$$

$$(T\Phi)(x) := \Phi_0(x) - \sum_{j=0}^{l} a_j(x)(I_{a+}^{\alpha-\alpha_j}\Phi)(x) . \tag{3.26}$$

Then the integral equation (3.8) can be rewritten as
$$\Phi(x) = (T\Phi)(x) . \tag{3.27}$$

Let

$$A = \max_{\substack{a \le x \le b \\ j \in {}_0Z_l}} |a_j(x)| . \tag{3.28}$$

Choose a positive number $\delta$ such that

$$\omega = A \sum_{j=0}^{l} \frac{\delta^{\alpha-\alpha_j}}{|\Gamma(\alpha-\alpha_j+1)|} < 1 . \tag{3.29}$$

Then we can prove the existence of a unique solution to the equation (3.8) in $L(a, x_1)$ with $x_1 = \min\{a+\delta, b\}$. By (3.26), we have

$$T : (a, x_1) \to L(a, x_1) \tag{3.30}$$

If $\Phi_1, \Phi_2 \in L(a, x_1)$, then from the boundedness of fractional integration operator we have

$$\|T\Phi_1 - T\Phi_2\|_{L(a, x_1)} \le \sum_{j=0}^{l} \left\| a_j(x)[I_{a+}^{\alpha-\alpha_j}(\Phi_1 - \Phi_2)](x) \right\|_{L(a, x_1)} \le$$

$$\le A \sum_{j=0}^{l} \frac{(x_1-a)^{\alpha-\alpha_j}}{|\Gamma(\alpha-\alpha_j+1)|} \|\Phi_1 - \Phi_2\|_{L(a, x_1)} = \omega \|\Phi_1 - \Phi_2\|_{L(a, x_1)}$$

and $0 < \omega < 1$, there exists the unique solution $\Phi(x) \in L(a, x_1)$ to the integral equation (3.8) on the interval $[a, x_1]$. And the solution $\Phi(x)$ is the limit of the sequence $\{(T^m \Phi_0^*)(x)\}$, i.e.

$$\lim_{m \to \infty} \left\| T^m \Phi_0^* - \Phi * \right\|_{L(a, x_1)} = 0. \tag{3.31}$$

Here $\Phi_0^*$ is a function in $L(a, b)$. If for at least one $k \in \{1, \cdots, n\}$, $b_k \neq 0$ in the initial condition, then we can take $\Phi_0^*(x) = \Phi_0(x)$.

$(T^m \Phi_0^*)(x)$ is written by the recursion formulas

$$(T^m \Phi_0^*)(x) = y_0(x) - \sum_{j=0}^{l} a_j(x) I_{a+}^{\alpha - \alpha_j} (T^{m-1} \Phi_0^*)(x). \tag{3.32}$$

If we denote $\Phi_m(x) := (T^m \Phi_0^*)(x)$, then we have

$$\Phi_m(x) = \Phi_0(x) - \sum_{j=0}^{l} a_j(x)(I_{a+}^{\alpha - \alpha_j} \Phi_{m-1})(x) \quad (m = 1, 2, \cdots). \tag{3.33}$$

(3.31) can be written as follows:

$$\lim_{m \to \infty} \left\| \Phi_m - \Phi \right\|_{L(a, x_1)} = 0. \tag{3.34}$$

If $x_1 \geq b$ then the proof is completed.

If $x_1 < b$ then let

$$\Phi_{01}(x) := \Phi_0(x) - \sum_{j=0}^{l} a_j(x)(I_{a+}^{\alpha - \alpha_j} \Phi)(x_1). \tag{3.35}$$

Note that $I_{a+}^{\alpha - \alpha_j} \Phi(x) = I_{a+}^{\alpha - \alpha_j} \Phi(x_1) + I_{x_1+}^{\alpha - \alpha_j} \Phi)(x)$. Then (3.8) is rewritten as follows

$$\Phi(x) = \Phi_{01}(x) - \sum_{j=0}^{l} a_j(x)(I_{x_1+}^{\alpha - \alpha_j} \Phi)(x). \tag{3.36}$$

From (3.29) and the above consideration, we can prove that the equation (3.36) has the unique solution $\Phi^* \in L(x_1, x_2)$ with $x_2 = \min\{x_1 + \delta, b\}$.

Repeating these processes, we conclude that the equation (3.8) has a unique solution $\Phi \in L(a, b)$.

Thus, there exists a unique solution $y(x) \in L^{\alpha}(a, b)$ to the Cauchy type problem (3.3),(3.6) by Theorem 3.1. By (3.16) we have

$$y(x) = \sum_{k=1}^{n} b_k H(k_0 - k) \frac{(x-a)^{\alpha - k}}{\Gamma(\alpha - k + 1)} + (I_{a+}^{\alpha} \Phi^*)(x) \tag{3.37}$$

This completes the proof of Theorem 3.2.

### 3.2  A Note for general Cauchy type problems

We proved that the special Cauchy type problem (3.3), (3.6) has the unique solution $y(x) \in L^{\alpha}(a, b)$. But in the case when the natural number $k_0$ defined by (3.5) is less than $n = [\text{Re}\,\alpha] + 1$, the Cauchy type problem

$$(D_{a+}^{\alpha} y)(x) + \sum_{j=0}^{l} a_j(x)(D_{a+}^{\alpha_j} y)(x) = g(x) \tag{3.38}$$

$$(D_{a+}^{\alpha - k} y)(a+) = b_k, \quad k = 1, \cdots, n \tag{3.39}$$

might not have solution $y(x) \in L^\alpha(a, b)$.

Under the assumption of Theorem 3.1, if we assume that there exists the solution $y(x) \in L^\alpha(a, b)$ of the problem (3.38), (3.38), then we have $y(x), (D_{a+}^\alpha y)(x) \in L(a, b)$ from the definition of $L^\alpha(a, b)$. The equation (3.38) can be rewritten as

$$(D_{a+}^\alpha y)(x) = g(x) - \sum_{j=0}^{l} a_j(x)(D_{a+}^{\alpha_j} y)(x). \tag{3.40}$$

Since $(D_{a+}^\alpha y)(x), g(x) \in L(a, b)$, we have

$$\sum_{j=0}^{l} a_j(x)(D_{a+}^{\alpha_j} y)(x) \in L(a, b). \tag{3.41}$$

By the definition of fractional derivatives, $(D_{a+}^\alpha y)(x) = (D^n I_{a+}^{n-\alpha} y)(x) \in L(a, b)$ and thus from (2.6) we have $y_{n-\alpha}(x) = (I_{a+}^{n-\alpha} y)(x) \in AC^n[a, b]$. From $n-1 < \operatorname{Re} \alpha < n$ and lemma 2.8 (a) we have

$$(I_{a+}^\alpha D_{a+}^\alpha y)(x) = y(x) - \sum_{k=1}^{n} y_{n-\alpha}^{(n-k)}(a+) \frac{(x-a)^{\alpha-k}}{\Gamma(\alpha-k+1)}.$$

From $n - k - 1 < \operatorname{Re}(\alpha - k) = \operatorname{Re}\alpha - k \leq n - k$ and the definition of fractional derivatives we have

$$y_{n-\alpha}^{(n-k)}(x) = (D^{n-k} I_{a+}^{n-\alpha} y)(x) = (D^{n-k} I_{a+}^{n-k-(\alpha-k)} y)(x) = (D_{a+}^{\alpha-k} y)(x).$$

Taking the operator $I_{a+}^\alpha$ to the both sides of (3.40) and considering (3.39), we have

$$y(x) = \sum_{k=1}^{n} b_k \frac{(x-a)^{\alpha-k}}{\Gamma(\alpha-k+1)} + I_{a+}^\alpha \left[ g(x) - \sum_{j=0}^{l} a_j(x)(D_{a+}^{\alpha_j} y)(x) \right]. \tag{3.42}$$

Operating the fractional derivative operator $D_{a+}^{\alpha_j}$ in the both sides and considering (2.20) and (2.28), we have

$$(D_{a+}^{\alpha_j} y)(x) = \sum_{k=1}^{n} b_k \frac{(x-a)^{\alpha-\alpha_j-k}}{\Gamma(\alpha-\alpha_j-k+1)} + I_{a+}^{\alpha-\alpha_j} \left[ g(x) - \sum_{j=0}^{l} a_j(x)(D_{a+}^{\alpha_j} y)(x) \right]$$

Therefore, we have

$$\sum_{j=0}^{l} a_j(x)(D_{a+}^{\alpha_j} y)(x) =$$

$$= \sum_{j=0}^{l} \sum_{k=1}^{k_0} b_k a_j(x) \frac{(x-a)^{\alpha-\alpha_j-k}}{\Gamma(\alpha-\alpha_j-k+1)} + \sum_{j=0}^{l} \sum_{k=k_0+1}^{n} b_k a_j(x) \frac{(x-a)^{\alpha-\alpha_j-k}}{\Gamma(\alpha-\alpha_j-k+1)} + \tag{3.43}$$

$$+ \sum_{j=0}^{l} a_j(x) I_{a+}^{\alpha-\alpha_j} \left[ g(x) - \sum_{j=0}^{l} a_j(x)(D_{a+}^{\alpha_j} y)(x) \right].$$

By the definition of $k_0$, the continuity of $a_j(x)$ ($j = 0, 1, \cdots, l$), the integrability of $g(x)$ and (3.41), the left side and the first and third terms of right side of (3.43) are integrable on $(a, b)$. If some $b_k$, $k \in \{k_0 + 1, \cdots, n\}$ are not zero, then some terms are not integrable in the second term of right side of (3.43). Because of these terms, in this case there might not exist the solution of

(3.38) and (3.39).

**Example 1.** For $b_i \in C \ (i = 1, 2, 3, 4)$, consider the problem

$$(D_{a+}^{3.5} y)(x) - 3(D_{0+}^{3.4} y)(x) = 0, \quad 0 < x < T \tag{3.44}$$

$$(D_{0+}^{2.5} y)(0+) = b_1, \ (D_{0+}^{1.5} y)(0+) = b_2, \ (D_{0+}^{0.5} y)(0+) = b_3, \ (I_{0+}^{0.5} y)(0+) = b_4. \tag{3.45}$$

If we assume that the problem (3.44) and (3.45) has a solution $y(x) \in L^{3.5}(0, T)$, then from (3.43) we have

$$(D_{0+}^{3.4} y)(x) = \frac{b_1 (x-0)^{-0.9}}{\Gamma(-1.9)} + \frac{b_2 (x-0)^{-1.9}}{\Gamma(-2.9)} + \frac{b_3 (x-0)^{-2.9}}{\Gamma(-3.9)} + \frac{b_4 (x-0)^{-3.9}}{\Gamma(-4.9)} + I_{0+}^{0.1}\left(D_{0+}^{3.4} y\right)(x).$$

Here all terms except for the second, third and fourth terms of right side are integrable and thus if the above equality holds, the the sum of these second, third and fourth terms of right side is integrable. This is only possible when $b_k$, $k \in \{2, 3, 4\}$ are all zero. So the problem of (3.44) and (3.45) has no solution if and only if some of $b_k$, $k \in \{2, 3, 4\}$ are not zero.

This example show us that if $a_j(x) \equiv a_j(\text{const})$, then the problem (3.38) and (3.39) with constant coefficients has a solution if and only if $b_k = 0$, $k \in \{k_0 + 1, \cdots n\}$.

### 3.3  Solving Method of the Cauchy type problems

In section 3.1, we saw that there exists a unique solution $y(x) \in L^\alpha(a, b)$ to the Cauchy type problem (3.3), (3.6), and the solution $y(x)$ is represented by (3.16) using the solution $\Phi(x) \in L(a, b)$ of the integral equation (3.8), where $\Phi(x) \in L(a, b)$ is obtained by successive approximations:

$$\Phi_0(x) = g(x) - \sum_{j=0}^{n} \sum_{k=1}^{n} b_k H(k_0 - k) a_j(x) \frac{(x-a)^{\alpha - \alpha_j - k}}{\Gamma(\alpha - \alpha_j - k + 1)},$$

$$\Phi_m(x) = \Phi_0(x) - \sum_{j=0}^{l} a_j(x)(I_{a+}^{\alpha - \alpha_j} \Phi_{m-1})(x), \quad m = 1, 2, \cdots, \tag{3.46}$$

$$\Phi(x) = \lim_{m \to \infty} \Phi_m(x) \ (\text{in } L(a, b)).$$

In this section, using this successive approximations, we get representation of the solution to the problem (3.3),(3.6) by coefficients and initial value of the equation.

### 3.3.1  Representation of the solution to the homogeneous equation.

We find a representation of the solution $y(x) \in L^\alpha(a, b)$ satisfying the linear homogeneous equation

$$(D_{a+}^\alpha y)(x) + \sum_{j=0}^{l} a_j(x)(D_{a+}^{\alpha_j} y)(x) = 0, \ x > a \tag{3.47}$$

and inhomogeneous initial condition (3.6).

**Definition 3.1.** Let $y_i(x) \in L^{\alpha}(a, b)$, $i = 1, \cdots, k_0$ be the solution of the problem

$$(D_{a+}^{\alpha} y_i)(x) + \sum_{j=0}^{l} a_j(x)(D_{a+}^{\alpha_j} y_i)(x) = 0, \ x > a, \tag{3.48}$$

$$(D_{a+}^{\alpha-k} y_i)(a+) = \delta_{ik}, \ k = 1, \cdots, n. \tag{3.49}$$

The system of $\{y_i(x) \in L^{\alpha}(a, b), i = 1, \cdots, k_0\}$ is called the canonical fundamental system of solutions of the linear homogeneous equation (3.47).

**Theorem 3.3.** Under the conditions of Theorem 3.1 for $\alpha$, $n$, $k_0$, $\alpha_j$, $a_j(x)$ ($j = 0, \cdots, l$), the canonical fundamental system $y_i(x) \in L^{\alpha}(a, b)$ ($i = 1, \cdots, k_0$) of solutions of the homogeneous equation (3.47) is given as follows:

$$y_i(x) = \frac{(x-a)^{\alpha-i}}{\Gamma(\alpha-i+1)} + \sum_{k=0}^{\infty}(-1)^{k+1} I_{a+}^{\alpha}\left(\sum_{j=0}^{l} a_j(x) I_{a+}^{\alpha-\alpha_j}\right)^k \sum_{j=0}^{l} a_j(x)\frac{(x-a)^{\alpha-\alpha_j-i}}{\Gamma(\alpha-\alpha_j-i+1)}. \tag{3.50}$$

**Proof.** By theorem 3.2, the Cauchy type problem (3.48),(3.49) has a unique solution $y_i(x) \in L^{\alpha}(a, b)$ ($i = 1, \cdots, k_0$). The integral equation (3.8) corresponding to the Cauchy type problem (3.48), (3.49) is

$$\Phi_i(x) = -\sum_{j=0}^{l} a_j(x)\frac{(x-a)^{\alpha-\alpha_j-i}}{\Gamma(\alpha-\alpha_j-i+1)} - \sum_{j=0}^{l} a_j(x)(I_{a+}^{\alpha-\alpha_j} \Phi_i)(x). \tag{3.51}$$

Form (3.16), by the solution $\Phi_i(x) \in L(a, b)$ of the equation (3.51), the solution $y_i(x) \in L^{\alpha}(a, b)$ of the problem (3.48) and (3.49) is represented by

$$y_i(x) = \frac{(x-a)^{\alpha-i}}{\Gamma(\alpha-i+1)} + (I_{a+}^{\alpha} \Phi_i)(x). \tag{3.52}$$

If we apply the successive approximations (3.46) to solve the solution of the equation (3.51), then we have

$$\Phi_{i0}(x) = -\sum_{j=0}^{l} a_j(x)\frac{(x-a)^{\alpha-\alpha_j-i}}{\Gamma(\alpha-\alpha_j-i+1)}$$

$$\Phi_{im}(x) = \Phi_{i0}(x) - \sum_{j=0}^{l} a_j(x)(I_{a+}^{\alpha-\alpha_j} \Phi_{im-1})(x), \ m = 1, 2, \cdots, \tag{3.53}$$

$$\Phi_i(x) = \lim_{m \to \infty} \Phi_{im}(x), \text{ in } L(a, b).$$

Calculating $\Phi_{i1}(x)$, we have

$$\Phi_{i1}(x) = \Phi_{i0}(x) - \sum_{j=0}^{l} a_j(x)(I_{a+}^{\alpha-\alpha_j} \Phi_{i0})(x) =$$

$$= -\sum_{j=0}^{l} a_j(x)\frac{(x-a)^{\alpha-\alpha_j-i}}{\Gamma(\alpha-\alpha_j-i+1)} + (-1)^2 \left(\sum_{j=0}^{l} a_j(x) I_{a+}^{\alpha-\alpha_j}\right)^1 \sum_{j=0}^{l} a_j(x)\frac{(x-a)^{\alpha-\alpha_j-i}}{\Gamma(\alpha-\alpha_j-i+1)}.$$

Calculating $\Phi_{i2}(x)$, we have

$$\Phi_{i2}(x) = \Phi_{i0} - \sum_{j=0}^{l} a_j(x)(I_{a+}^{\alpha-\alpha_j}\Phi_{i1})(x) = \Phi_{i0}(x) - \sum_{j=0}^{l} a_j(x) I_{a+}^{\alpha-\alpha_j}$$

$$\left( -\sum_{j=0}^{l} a_j(x) \frac{(x-a)^{\alpha-\alpha_j-i}}{\Gamma(\alpha-\alpha_j-i+1)} + (-1)^2 \left( \sum_{j=0}^{l} a_j(x) I_{a+}^{\alpha-\alpha_j} \right)^1 \sum_{j=0}^{l} a_j(x) \frac{(x-a)^{\alpha-\alpha_j-i}}{\Gamma(\alpha-\alpha_j-i+1)} \right) =$$

$$= -\sum_{j=0}^{l} a_j(x) \frac{(x-a)^{\alpha-\alpha_j-i}}{\Gamma(\alpha-\alpha_j-i+1)} + (-1)^2 \left( \sum_{j=0}^{l} a_j(x) I_{a+}^{\alpha-\alpha_j} \right)^1 \sum_{j=0}^{l} a_j(x) \frac{(x-a)^{\alpha-\alpha_j-i}}{\Gamma(\alpha-\alpha_j-i+1)} +$$

$$+ (-1)^3 \left( \sum_{j=0}^{l} a_j(x) I_{a+}^{\alpha-\alpha_j} \right)^2 \sum_{j=0}^{l} a_j(x) \frac{(x-a)^{\alpha-\alpha_j-i}}{\Gamma(\alpha-\alpha_j-i+1)} =$$

$$= \sum_{k=0}^{2} (-1)^{k+1} \left( \sum_{j=0}^{l} a_j(x) I_{a+}^{\alpha-\alpha_j} \right)^k \sum_{j=0}^{l} a_j(x) \frac{(x-a)^{\alpha-\alpha_j-i}}{\Gamma(\alpha-\alpha_j-i+1)}$$

Calculating $\Phi_{im}(x)$ similarly, we have

$$\Phi_{im}(x) = \sum_{k=0}^{m} (-1)^{k+1} \left( \sum_{j=0}^{l} a_j(x) I_{a+}^{\alpha-\alpha_j} \right)^k \sum_{j=0}^{l} a_j(x) \frac{(x-a)^{\alpha-\alpha_j-i}}{\Gamma(\alpha-\alpha_j-i+1)}.$$

Therefore we have

$$\Phi_i(x) = \lim_{\substack{m \to \infty \\ L(a,b)}} \Phi_{im}(x) = \sum_{k=0}^{\infty} (-1)^{k+1} \left( \sum_{j=0}^{l} a_j(x) I_{a+}^{\alpha-\alpha_j} \right)^k \sum_{j=0}^{l} a_j(x) \frac{(x-a)^{\alpha-\alpha_j-i}}{\Gamma(\alpha-\alpha_j-i+1)}. \quad (3.54)$$

Here $\left( \sum_{j=0}^{l} a_j(x) I_{a+}^{\alpha-\alpha_j} \right)^k$ is the composition operator of $k$ times of the operator $\sum_{j=0}^{l} a_j(x) I_{a+}^{\alpha-\alpha_j}$ and in case with $k=0$ is a unit operator. Substituting (3.54) into (3.52), we have (3.50). The proof of the theorem is completed.

**Corollary 3.4.** Under the conditions of Theorem 3.1 for $\alpha$, $n$, $k_0$, $\alpha_j, a_j(x)$ $(j=0,\cdots,l)$, the solution $y(x) \in L^\alpha(a, b)$ to the Cauchy type problem (3.3), (3.6) is represented as follows:

$$y(x) = \sum_{i=1}^{k_0} b_i y_i(x). \quad (3.55)$$

Here $y_i(x)$ $(i=1,\cdots,k_0)$ is the canonical fundamental system of the solution to the homogeneous equation (3.47).

**Remark 3.1**. (3.55) with complex numbers $b_i \in C$ $(i=1,\cdots,k_0)$ is the general solution to the homogeneous equation (3.3).

### 3.3.2 Representation of the solution to the inhomogeneous equation

We find a representation of the solution $y(x) \in L^\alpha(a, b)$ to the problem

$$(D_{a+}^{\alpha} y)(x) + \sum_{j=0}^{l} a_j(x)(D_{a+}^{\alpha_j} y)(x) = g(x), \ x > a \tag{3.56}$$

$$(D_{a+}^{\alpha-k} y)(a+) = 0, \ k = 1, \cdots, n. \tag{3.57}$$

**Theorem 3.5.** Under the conditions of Theorem 3.1 for $\alpha$, $n$, $k_0$, $\alpha_j$, $a_j(x)$ ($j = 0, \cdots, l$), $g(x)$, the solution $y(x) \in L^{\alpha}(a, b)$ to the problem (3.56)-(3.57) is represented by

$$y(x) = \sum_{k=0}^{\infty} (-1)^k I_{a+}^{\alpha} \left( \sum_{j=0}^{l} a_j(x) I_{a+}^{\alpha-\alpha_j} \right)^k g(x). \tag{3.58}$$

**Proof.** By Theorem 3.2, this problem has a unique solution $y(x) \in L^{\alpha}(a, b)$. This solution is obtained by (3.8), (3.16) and (3.46). The integral equation (3.8) corresponding to the problem (3.56) and (3.57) is

$$\Phi(x) = g(x) - \sum_{j=0}^{l} a_j(x)(I_{a+}^{\alpha-\alpha_j} \Phi)(x). \tag{3.59}$$

The corresponding (3.16) is

$$y(x) = (I_{a+}^{\alpha} \Phi)(x). \tag{3.60}$$

The corresponding successive approximation is given by

$$\Phi_0(x) = g(x)$$

$$\Phi_m(x) = \Phi_0(x) - \sum_{j=0}^{l} a_j(x)(I_{a+}^{\alpha-\alpha_j} \Phi_{m-1})(x), \ m = 1, 2, \cdots \tag{3.61}$$

$$\Phi(x) = \lim_{m \to \infty} \Phi_m(x) \text{ in } L(a, b).$$

In the similar way with theorem 3.3 deriving of $y_i(x)$, $i = 1, \cdots, k_0$, (3.58) is derived by (3.59), (3.60) and (3.61). The proof of the theorem is completed.

**Definition 3.2.** The solution $G(x; \xi) \in L^{\alpha}(\xi, b)$, $\xi > a$ to the Cauchy type problem;

$$(D_{\xi+}^{\alpha} G)(x; \xi) + \sum_{j=0}^{l} a_j(x)(D_{\xi+}^{\alpha_j} G)(x; \xi) = 0, \ x > \xi \tag{3.62}$$

$$(D_{\xi+}^{\alpha-k} G)(x; \xi)\Big|_{x=\xi} = \begin{cases} 1, & k = 1 \\ 0, & k \neq 1 \end{cases}, \ k = 1, \cdots, n \tag{3.63}$$

is called the Green's function of the inhomogeneous Cauchy type problem (3.56),(3.57).

**Remark 3.2.** According to the definition, the Green's function is just the first function $y_1(x)$ of the canonical fundamental system of the solution to the homogeneous equation when $a = \xi$ and thus there exists the Green's function $G(x; \xi) \in L^{\alpha}(\xi, b)$ of the Cauchy type problem (3.56) (3.57) and it is represented as follows:

$$G(x; \xi) = \frac{(x-\xi)^{\alpha-1}}{\Gamma(\alpha)} + \sum_{k=0}^{\infty} (-1)^{k+1} I_{\xi+}^{\alpha} \left( \sum_{j=0}^{l} a_j(x) I_{\xi+}^{\alpha-\alpha_j} \right)^k \sum_{j=0}^{l} a_j(x) \frac{(x-\xi)^{\alpha-\alpha_j-1}}{\Gamma(\alpha-\alpha_j)}. \tag{3.64}$$

**Theorem 3.6.** Under the conditions of Theorem 3.1 for $\alpha$, $n$, $k_0$, $\alpha_j$, $a_j(x)$ ($j = 0, \cdots, l$), $g(x)$, the solution $y(x) \in L^{\alpha}(a, b)$ to the inhomogeneous Cauchy type problem (3.56) (3.57) is

represented as follows:

$$y(x) = \int_a^x G(x;\ \xi) g(\xi) d\xi. \tag{3.65}$$

And (3.65) is equal to the solution given by (3.58).

**Proof.** In Theorem 3.5 we proved that the unique solution $y(x) \in L^\alpha(a,\ b)$ to the inhomogeneous Cauchy type problem (3.56), (3.57) is represented by (3.58). Therefore, if we prove that (3.65) is equal to (3.58), the proof is completed.

By the representation of the Green's function $G(x;\ \xi)$,

$$y(x) = \int_a^x G(x;\ \xi) g(\xi) d\xi = \int_a^x \frac{(x-\xi)^{\alpha-1}}{\Gamma(\alpha)} g(\xi) d\xi +$$

$$+ \sum_{k=0}^n (-1)^{k+1} \int_a^x I_{\xi+}^\alpha \left( \sum_{j=0}^l a_j(x) I_{\xi+}^{\alpha-\alpha_j} \right)^k \sum_{j=0}^l a_j(x) \frac{(x-\xi)^{\alpha-\alpha_j-1}}{\Gamma(\alpha-\alpha_j)} g(\xi) d\xi.$$

By the definition of fractional integral operators $I_{\xi+}^\alpha$, $I_{\xi+}^{\alpha-\alpha_j}$ and the generalized Fubini's theorem, we have

$$y(x) = (I_{a+}^\alpha g)(x) + \sum_{k=0}^\infty (-1)^{k+1} I_{a+}^\alpha \left( \sum_{j=0}^l a_j(x) I_{a+}^{\alpha-\alpha_j} \right)^{k+1} g(x) =$$

$$= \sum_{k=0}^\infty (-1)^k I_{a+}^\alpha \left( \sum_{j=0}^l a_j(x) I_{a+}^{\alpha-\alpha_j} \right)^k g(x).$$

**Theorem 3.7.** Under the conditions of Theorem 3.1 for $\alpha$, $n$, $k_0$, $\alpha_j, a_j(x)$ $(j=0,\cdots,\ l)$, $g(x), b_k \in C$ $(k=1,\cdots,\ n)$, the unique solution $y(x) \in L^\alpha(a,\ b)$ to the problem (3.3) (3.6) is represented as follows

$$y(x) = \sum_{i=1}^{k_0} b_i y_i(x) + \int_a^x G(x;\ \xi) g(\xi) d\xi. \tag{3.66}$$

Here $\{\ y_i(x)\ (i=1,\cdots,\ k_0)\ \}$ is the canonical fundamental system (3.50) of the solutions to the homogeneous equation (3.47) and the function $G(x;\ \xi)$ is the Green's function (3.64) of the inhomogeneous Cauchy type problem (3.56) (3.57).

### 3.3.3 Note on linear fractional differential equations with constant coefficients

We consider the linear fractional differential equation with constant coefficients;

$$(D_{a+}^\alpha y)(x) + \sum_{j=0}^l a_j (D_{a+}^{\alpha_j} y)(x) = g(x),\ x > a \tag{3.67}$$

$$(D_{a+}^{\alpha-k} y)(a+) = b_k H(k_0 - k),\ k = 1,\ 2,\cdots,\ n. \tag{3.68}$$

Here $a_j$ $(j = 0,\ 1,\cdots,\ l)$ are the complex constants.

The homogeneous equation corresponding to the inhomogeneous equation (3.67) is

$$(D_{a+}^{\alpha} y)(x) + \sum_{j=0}^{l} a_j (D_{a+}^{\alpha_j} y)(x) = 0. \tag{3.69}$$

Now we consider the representation of the canonical fundamental system of the solutions to the homogeneous equation (3.69). As a special case of the results of Theorem 3.3 we have the following theorem.

**Theorem 3.8.** Under the conditions of Theorem 3.1 for $\alpha \in C$, $\alpha_j \in C$ ($j = 0, 1, \cdots, l$), $k_0$, $n$, the canonical fundamental system $y_i(x) \in L^{\alpha}(a, b)$ ($i = 1, \cdots, k_0$) of the solutions to the homogeneous equation (3.69) is represented as follows:

$$y_i(x) = \frac{(x-a)^{\alpha-i}}{\Gamma(\alpha-i+1)} + \sum_{k=0}^{\infty}(-1)^{k+1} \sum_{|\beta|=k} \frac{k!}{\beta_0! \cdots \beta_l!} a_0^{\beta_0} \cdots a_l^{\beta_l} I_{a+}^{(\alpha-\alpha_0)\beta_0+\cdots+(\alpha-\alpha_l)\beta_l+\alpha}$$

$$\cdot \sum_{j=0}^{l} a_j \frac{(x-a)^{\alpha-\alpha_j-i}}{\Gamma(\alpha-\alpha_j-i+1)}, \quad i = 1, 2, \cdots, k_0. \tag{3.70}$$

Here $\beta = (\beta_0, \cdots, \beta_l) \in Z_+^{l+1}$ and $|\beta| = \beta_0 + \cdots + \beta_l$.

**Remark 3.3.** The Green's function $G(x; \xi) \in L^{\alpha}(\xi, b)$ for the Cauchy type problem (3.67), (3.68) is represented as follows:

$$G(x; \xi) = \sum_{k=0}^{\infty}(-1)^k \sum_{|\beta|=k} \frac{k!}{\beta_0! \cdots \beta_l!} a_0^{\beta_0} \cdots a_l^{\beta_l} \frac{(x-a)^{(\alpha-\alpha_0)\beta_0+\cdots+(\alpha-\alpha_l)\beta_l+\alpha-1}}{\Gamma((\alpha-\alpha_0)\beta_0+\cdots+(\alpha-\alpha_l)\beta_l+\alpha)}. \tag{3.71}$$

**Remark 3.4.** If we let $\xi \to 0+$ in (3.71), then we have

$$G(x) = \lim_{\xi \to 0+} G(x; \xi) = \sum_{k=0}^{\infty}(-1)^k \sum_{|\beta|=k} \frac{k!}{\beta_0! \cdots \beta_l!} a_0^{\beta_0} \cdots a_l^{\beta_l} \frac{x^{(\alpha-\alpha_0)\beta_0+\cdots+(\alpha-\alpha_l)\beta_l+\alpha-1}}{\Gamma((\alpha-\alpha_0)\beta_0+\cdots+(\alpha-\alpha_l)\beta_l+\alpha)}. \tag{3.72}$$

This is equal to the result of [29] ( See the equation (2.25) of [29] )

**Remark 3.5.** Under the conditions of Theorem 3.1 for $\alpha \in C$, $\alpha_j \in C$ ($j = 0, 1, \cdots, l$), $n \in N$, $g(x)$, $b_k \in C$ ($k = 1, \cdots, n$), $k_0$, the unique solution $y(x) \in L^{\alpha}(a, b)$ to the problem (3.67)-(3.68) is represented as follows:

$$y(x) = \sum_{i=1}^{k_0} b_i y_i(x) + \int_a^x G(x; \xi) g(\xi) d\xi. \tag{3.73}$$

Here $y_i(x)$ ($i = 1, \cdots, k_0$) and $G(x; \xi)$ are equal to (3.70) and (3.71), respectively. This result is obtained in [29] using Neuman's series method in the case of $a = 0$. (See (2.22) and theorem 1 in [29].)

### 3.3.4 Some examples.

**Example 2.** Consider the following linear homogeneous fractional differential equation with constant coefficients

$$(D_{0+}^{1.5} y)(x) + 3(D_{0+}^1 y)(x) = 0, \quad 0 < x < T, \quad T > 0 \tag{3.74}$$

This is the case when $\alpha = 1.5$, $n = 2$, $\alpha_1 = 1$, $l = 1$, $a_1 = 3$, $a = 0$, $b = T$ in (3.3). We can obtain the $k_0$ by the equation (3.5). Since

$$a_1 \frac{x^{\alpha-\alpha_1-k}}{\Gamma(\alpha-\alpha_1-k+1)} = 3\frac{x^{0.5-k}}{\Gamma(0.5-k+1)} \begin{cases} \in L(0,\ T),\ k=1, \\ \notin L(0,\ T),\ k=2=n, \end{cases}$$

we have $k_0 = 1 < 2 = n$. There exists a unique solution $y(x) \in L^{1.5}(0,\ T)$ to the equation (3.74) satisfying the following initial conditions

$$(D_{0+}^{1.5-1} y)(0+) = (D_{0+}^{0.5} y)(0+) = b_1,$$
$$(D_{0+}^{1.5-2} y)(0+) = (D_{0+}^{-0.5} y)(0+) = 0. \tag{3.75}$$

Then the canonical fundamental system $y_1(x)$ of the solution is provided by

$$y_1(x) = \sum_{k=0}^{\infty} (-1)^k 3^k \frac{x^{1.5k+0.5}}{\Gamma(1.5k+1.5)}$$

using the formula (3.70), and the solution is represent as $y(x) = b_1 y_1(x)$. But there doesn't exist the solution $y(x) \in L^{1.5}(0,\ T)$ to the equation (3.74) satisfying the initial conditions

$$(D_{0+}^{0.5} y)(0+) = b_1,\quad (D_{0+}^{-0.5} y)(0+) = b_2 \tag{3.76}$$

with any $b_1, b_2 \in C$ and $b_2 \neq 0$.

**Example 3.** Consider the following linear homogeneous fractional differential equation with variable coefficients

$$(D_{0+}^{1.5} y)(x) + x(D_{0+}^{1} y)(x) = 0,\ 0 < x < T. \tag{3.77}$$

This is the case when $\alpha = 1.5$, $n = 2$, $\alpha_1 = 1$, $l = 1$, $a_1(x) = x$, $a = 0$, $b = T$ in (3.3). We can obtain the $k_0$ by the equation (3.5). Since

$$a_1(x) \frac{x^{\alpha-\alpha_1-k}}{\Gamma(\alpha-\alpha_1-k+1)} = x\frac{x^{1.5-1-k}}{\Gamma(1.5-1-k+1)} = \frac{x^{1.5-k}}{\Gamma(1.5-k)} \in L(0,\ T),\ k=1,\ 2,$$

we have $k_0 = 2 = n$. By the formula (3.50), the canonical fundamental system $y_i(x) \in L^{1.5}(0,\ T)$, $i = 1,\ 2$ of the solutions to (3.77) is represented as follows:

$$y_i(x) = \frac{x^{1.5-i}}{\Gamma(1.5-i+1)} + \sum_{k=0}^{\infty} (-1)^{k+1} I_{0+}^{1.5}(xI_{0+}^{1.5-1})^k x\frac{x^{1.5-1-i}}{\Gamma(1.5-1-i+1)} =$$
$$= \frac{x^{1.5-i}}{\Gamma(2.5-i)} + \sum_{k=0}^{\infty} (-1)^{k+1} \prod_{j=1}^{k+1}(1.5j-i)\frac{x^{(k+2)1.5-i}}{\Gamma((k+2)1.5-i+1)},\ i=1,\ 2. \tag{3.78}$$

The unique solution $y(x) \in L^{1.5}(0,\ T)$ to (3.77) satisfying the initial conditions

$$(D_{0+}^{0.5} y)(0+) = b_1,\quad (D_{0+}^{-0.5} y)(0+) = b_2 \tag{3.79}$$

for any $b_i \in C$, $i = 1,\ 2$ is represented as follows:

$$y(x) = b_1 y_1(x) + b_2 y_2(x). \tag{3.80}$$

**Example 4.** Consider the following linear homogeneous fractional differential equation with constant coefficients

$$(D_{0+}^{3.5} y)(x) + 2(D_{0+}^{2.5} y)(x) + 3(D_{0+}^{1.5} y)(x) + 5(D_{0+}^{1.3} y)(x) + 4(D_{0+}^{1.2} y)(x) +$$
$$+ 6(D_{0+}^{0.5} y)(x) + 9(D_{0+}^{0.4} y)(x) + 7y(x) = 0. \tag{3.81}$$

This is the case when

$\alpha = 3.5$, $n = 4$, $l = 6$, $a_6 = 2$, $a_5 = 3$, $a_4 = 5$, $a_3 = 4$, $a_2 = 6$, $a_1 = 9$, $a_0 = 7$;
$\alpha_6 = 2.5$, $\alpha_5 = 1.5$, $\alpha_4 = 1.3$, $\alpha_3 = 1.2$, $\alpha_2 = 0.5$, $\alpha_1 = 0.4$, $\alpha_0 = 0$

in (3.3). We use the following table to decide $k_0$ using (3.5).

| $k$ | $a_6 = 2$ | $a_5 = 3$ | $a_4 = 5$ | $a_3 = 4$ | $a_2 = 6$ | $a_1 = 9$ | $a_0 = 7$ | $a_j \dfrac{x^{\alpha - \alpha_j - k}}{\Gamma(\alpha - \alpha_j - k + 1)}$ |
|---|---|---|---|---|---|---|---|---|
|   | $\alpha - \alpha_6 - k + 1$ | $\alpha - \alpha_5 - k + 1$ | $\alpha - \alpha_4 - k + 1$ | $\alpha - \alpha_3 - k + 1$ | $\alpha - \alpha_2 - k + 1$ | $\alpha - \alpha_1 - k + 1$ | $\alpha - \alpha_0 - k + 1$ | $j = 0, 1, \cdots, 6$ |
| 1 | 1 | 2 | 2.2 | 2.3 | 3 | 3.1 | 3.5 | $\in L(0, T)$ |
| 2 | 0 | 1 | 1.2 | 1.3 | 2 | 2.1 | 2.5 | $\in L(0, T)$ |
| 3 | -1 | 0 | 0.2 | 0.3 | 1 | 1.1 | 1.5 | $\in L(0, T)$ |
| 4 | -2 | -1 | (-0.8) | (-0.7) | 0 | 0.1 | 0.5 | $\notin L(0, T)$ |

By this table, $k_0 = 3 < 4 = n$. By the formula (3.70), the canonical fundamental system of the solutions is provided as follows

$$y_1(x) = \frac{x^{2.5}}{\Gamma(3.5)} + \sum_{k=0}^{\infty} (-1)^{k+1} \sum_{|\beta|=k} \frac{k!}{\beta_0! \beta_1! \cdots \beta_6!} 7^{\beta_0} 9^{\beta_1} 6^{\beta_2} 4^{\beta_3} 5^{\beta_4} 3^{\beta_5} 2^{\beta_6} \cdot$$
$$\cdot I_{0+}^{3.5\beta_0 + 3.1\beta_1 + 3\beta_2 + 2.3\beta_3 + 2.2\beta_4 + 2\beta_5 + 1\beta_6 + 3.5}$$
$$\left( 7 \frac{x^{2.5}}{\Gamma(3.5)} + 9 \frac{x^{2.1}}{\Gamma(3.1)} + 6 \frac{x^2}{\Gamma(3)} + 4 \frac{x^{1.3}}{\Gamma(2.3)} + 5 \frac{x^{1.2}}{\Gamma(2.2)} + 3 \frac{x}{\Gamma(2)} + 2 \frac{x^0}{\Gamma(1)} \right).$$

$$y_2(x) = \frac{x^{1.5}}{\Gamma(2.5)} + \sum_{k=0}^{\infty} (-1)^{k+1} \sum_{|\beta|=k} \frac{k!}{\beta_0! \beta_1! \cdots \beta_6!} 7^{\beta_0} 9^{\beta_1} 6^{\beta_2} 4^{\beta_3} 5^{\beta_4} 3^{\beta_5} 2^{\beta_6} \cdot$$
$$\cdot I_{0+}^{3.5\beta_0 + 3.1\beta_1 + 3\beta_2 + 2.3\beta_3 + 2.2\beta_4 + 2\beta_5 + 1\beta_6 + 3.5}$$
$$\left( 7 \frac{x^{1.5}}{\Gamma(2.5)} + 9 \frac{x^{1.1}}{\Gamma(2.1)} + 6 \frac{x}{\Gamma(2)} + 4 \frac{x^{0.3}}{\Gamma(1.3)} + 5 \frac{x^{0.2}}{\Gamma(1.2)} + 3 \frac{x^0}{\Gamma(1)} \right).$$

$$y_3(x) = \frac{x^{0.5}}{\Gamma(1.5)} + \sum_{k=0}^{\infty} (-1)^{k+1} \sum_{|\beta|=k} \frac{k!}{\beta_0! \beta_1! \cdots \beta_6!} 7^{\beta_0} 9^{\beta_1} 6^{\beta_2} 4^{\beta_3} 5^{\beta_4} 3^{\beta_5} 2^{\beta_6} \cdot$$
$$\cdot I_{0+}^{3.5\beta_0 + 3.1\beta_1 + 3\beta_2 + 2.3\beta_3 + 2.2\beta_4 + 2\beta_5 + 1\beta_6 + 3.5}$$
$$\left( 7 \frac{x^{0.5}}{\Gamma(1.5)} + 9 \frac{x^{0.1}}{\Gamma(1.1)} + 6 \frac{x^0}{\Gamma(1)} + 4 \frac{x^{-0.7}}{\Gamma(0.3)} + 5 \frac{x^{-0.8}}{\Gamma(0.2)} \right)$$

Then the unique solution $y(x) \in L^{3.5}(0, T)$ to (3.81) satisfying the initial conditions

$$(D_{0+}^{3.5-k} y)(0+) = b_k, \quad k = 1, \cdots, 4 \qquad \qquad 3.82)$$

with $b_i \in C$ ($i = 1, 2, 3, 4$) and $b_4 = 0$ is represented as follows:

$$y(x) = b_1 y_1(x) + b_2 y_2(x) + b_3 y_3(x).$$

But there doesn't exist the solution $y(x) \in L^{3.5}(0, T)$ to (3.81) satisfying the initial conditions

$$(D_{0+}^{3.5-k} y)(0+) = b_k, \; k = 1, \cdots, 4 \tag{3.83}$$

for $b_i \in C \; (i = 1, 2, 3, 4)$ and $b_4 \neq 0$.

**Example 3.5** Consider the following linear homogeneous fractional differential equation of complex order with constant coefficients

$$(D_{a+}^{3.5+i2.6} y)(x) - 3(D_{0+}^{3.4+i2.6} y)(x) = 0, \; 0 < x < T. \tag{3.84}$$

In this case $k_0 = 1$ and the canonical fundamental system of the solutions is represented by

$$y_1(x) = \sum_{k=0}^{\infty} 3^k \frac{x^{0.1k+2.5+i2.6}}{\Gamma(0.1k + 3.5 + i2.6)}. \tag{3.85}$$

The unique solution $y(x) \in L^{3.5}(0, T)$ to (3.84) satisfying the following initial condition

$$(D_{0+}^{2.5} y)(0+) = b_1, \; (D_{0+}^{1.5} y)(0+) = (D_{0+}^{0.5} y)(0+) = (D_{0+}^{-0.5} y)(0+) = 0. \tag{3.86}$$

with complex number $b_1 \in C$ is represented by $y(x) = b_1 y_1(x)$. But by Theorem 3.3 there doesn't exist the solution $y(x) \in L^{3.5}(0, T)$ to the equation (3.84) satisfying the initial conditions

$$(D_{0+}^{3.5-k} y)(0) = b_k, \; k = 1, \cdots, 4 \tag{3.87}$$

with $b_k \in C, \; k = 2, 3, 4$ such that at least one of them are not equal to zero.

## 4. Conclusion

Our examples 1, 2, 4 and 5 are counter examples of corollary 3.6 at page 158 of [16] which asserts that the problem (3.38) and (3.39) with all $b_k$, $k = 1, \cdots, n$ has a unique solution in $L^\alpha[a, b]$.

According to our results, the Cauchy type problem (3.38) and (3.39) might not have a solution in the case when some of $b_k \in C$, $k = k_0 + 1, \cdots, n$ are not zero.

We provided a solving of method for the Cauchy type problem with $b_k = 0$, $k = k_0 + 1, \cdots, n$ and the solution representations. Our method is just the same method of [16]. We only took notice that some terms are not integrable.

## References

[1] Bonilla,B.,Kilbas,A.A., and Trujillo,J.J., Systems of nonlinear fractional differential equations in the space of summable functions,Tr.Inst.Mat.Minsk,6(2000) 38-46
[2] Bonilla,B.,Rivero,M. and Trujillo,J.J.,On systems linear fractional differential equations with constant coefficients,Applied Mathematics and Computation,187(2007)68-78
[3] Caponetto,R.,Dongola,G.,Fortuna,L.,Petras, Fractional order systems:Modeling and Control Application,in : World Scientific Series on Nonlinear, Science,Vol.72,World Scientific,2010
[4] Fukunaga,M. and Shimizu,N., Analytical and numerical Solutions for Fractional Viscoelastic Equations, JSME Int. J., Series C,47(2004)251-259
[5] Garra, R., Analytic solution of a class of fractional differential equations with variable coefficients by operatorial methods, Commun Nonlinear Sci Numer Simulat 17 (2012) 1549–1554
[6] Gorenflo,R.,Luchko,Yu., Operational method for solving generalized Abel integral equations of second kind integral transforms and special functions,5(1997)47-58
[7] Hadid,S.B.,Luchko,Yu.F., An Operational method for solving fractional differential equations of an arbitrary real order,Panamer.Math.J.,6(1996)57-63

[8] Hilfer,R.,Applications of Fractional Calculus in Physics, World Scientific, Singapore(2000)

[9] Hilfer,H.,Luchko,Y. and Tomovski,Z., Operational method for the solution of fractional differential equations with generalized Riemann-Liouville fractional derivatives, Fract.Calc.Appl.Anal.,12,3(2009)299-318

[10] Hu,Y.,Junshong,D., Analytical Solution of the linear fractional differential equation by Adomian decomposition method,Comput. and appl. Math.,(2007)doi:10.1016/j cam.14.005

[11] Kilbas,A.A., Bonilla,B., and Trujillo,J.J., Existence and uniqueness theorems for nonlinear fractional differential equations , Demonstr.Math.,33(2.3),(2000)583-602

[12] Kilbas,A.A., Bonilla,B., and Trujillo,J.J.,Fractional integrals and derivatives , and differential equations of fractional order in weighted spaces of continuous functions(Russian), Dokl.Nats.Akad.Nauk Belarusi, 44(2.6),(2000) 18-22

[13] Kilbas,A.A., Bonilla,B., and Trujillo,J.J., Nonlinear differential equations of fractional order in space of integrable functions(Russian),Dokl.Akad.Nauk,374(2.4),(2000),445-449,Translated in Dokl. Math.,62(2.2) (2000),222-226

[14] Kilbas,A.A.,Rivero,M.,Rodrignez-Germa L.,Trujillo J.J., $\alpha$-Analytic solutions of some linear fractional differential equations with variable coefficients,Appl. Math. and Comput.,187(2007)239-249

[15] Kilbas,A.A., Saigo,M., On Mittag-Leffler type functions,Fractional Calculus Operations and Solution of integral equations,Integral Transforms and Special Functions,4(1996)355-370

[16] Kilbas,A.A.,Srivastava,H.M. and Trujillo,J.J., Theory and Applications of Fractional Differential Equations,Elsevier,Amsterdam-Tokyo,(2006)

[17] V. Kiryakova, Generalized Fractional Calculus and Applications. Long-man & J. Wiley, Harlow & N. York (1994).

[18] Luchko,Y., Operational method in fractional calculus, Fract.Calc.Appl.Anal., 2,4(1999)463-468

[19] Luchko,Y. and Gorenflo,R., An operational method for solving fractional differential equations,Acta Mathematica Vietnamica 24(1999),207-234

[20] Luchko,Yu.F. and Srivastava,H.M., The Exact Solution of Certain Differential Equations of Fractional Order by Using Operational Calculus,Computers Math.Applic.,29(1995)73-85

[21] Machado,J.J.,Kiryakova,V.,Mainardi,F., Recent history of fractional calculus, Commun.Nonlinear Sci. Numer. Simul.,16(2.3)(2011),1140-1153

[22] Mainardi,F.,Fractional relaxation and fractional diffusion equations,mathematical aspects,In:Proceedings of the 12-th IMACS World Congress(Ed.W.F.Ames),Georgia Tech Atlanta,Vol.1,(1994),329-332

[23] Miller,K.S. and Ross,B., An Introduction to the Fractional Calculus and Fractional Differential Equations,Wiley and Sons,New York,1993

[24] Momani,S. and Obidat,Z., Numerical comparison of methods for solving linear differential equations of fractional order,Chaos,Solitons and Fractals,31(2007)1248-1255

[25] Morita,T., Solution of the Initial-Value Problem of a Fractional Differential Equation,int.J.Appl.Math.,15(2004)299-313

[26] Morita,T.,"Solution of the Initial-Value problem of a Fractional Differential Equation",in Similarity in Diversity,Eds.S.Fujita,H.Hara,D.Morabito and Y.Okamura,Nova Science Publishes,Inc.,New York(2003)245-257

[27] Morita,T. and Sato,K.,Neumann-Series Solution of Fractional Differential Equation, Interdisciplinary Information Sciences,16,1(2010) 127-137

[28] Morita,T. and Sato,K., Solution of Fractional Differential Equation in Terons of Distribution Theory, Interdisc. Inf. Sc., 12,2(2006) 71-83

[29] Myong-Ha Kim, Guk-Chol Ri and Hyong-Chol O, Operation Method for Solving Multi-term fractional Differential Equations with the Generalized fractional Derivatives, Fract. Calc. Appl. Anal., 17, 1(2014) 79-95

[30] Myong-Ha Kim, Hyong-Chol O, Explicit Representation of Green's Function for Linear Fractional Differential operator with variable Coefficients.,Journal of fractional Calculus and Applicationx, 5,1(2014), 26-28

[31] Podlubny,I., Fractional Differential Equations,Academic Press,San Diego,1999

[32] Podlubny,I., Solution of linear fractional differential equations with constant coefficients,In:Transform Methods and Special Functions (EDS, P.Rusev, I.Dimovski, V. Kiryakova), Singapore,Science Culture Technology Publishing(SCTP)(1995)227-237

[33] Podlubny,I., The Laplace Transform Method for Linear Differential Equations of the Fractional order,Slovak Academy of Sciences,Institute of Experimental Physics, arxiv:funct-am/97|V005v.I 30Oct.1997

[34] Podlubny,I., The Laplace Transform Method for Linear Differential Equations of the Fractional order, Inst.Exp.Phys.,Slovak Acad.Sci.No UEF-02-94,1994,Kosice,1-32

[35] Samko,S.G.,Kilbas,A.A.,Marichev,O.I., Fractional Intgrals and Derivatives:Theory and Applications,New York and London,Gordon and Breach Science Publishers(1993)

[36] Sandev,T.,Metzler,R.,Tomovski,Z., Fractional diffusion equation with a generalized Riemann-Liouville time fractional derivative, J. Phys. A:Math.Theor., 44(2.25)(2011) Article ID 255203,21p

[37] Tomovski,Z., Generalized Cauchy type problems for nonlinear fractional differential equations with composite fractional derivative operator,Nonlinear Analysis75(2012) 3364-3384